\documentclass[12pt]{article}
\usepackage{amssymb}
\usepackage{amsmath}
\usepackage{graphicx}
\usepackage{cite}
\usepackage{xcolor}
  \usepackage{paralist}
  \usepackage{graphics}
  \usepackage{epsfig} 
  \usepackage{bigints}
  \usepackage{bm}
  \usepackage{epstopdf}
\usepackage{mathtools}

\newtheorem{theorem}{Theorem}[section]

\newtheorem{lemma}[theorem]{Lemma}

\newtheorem{remark}{Remark}

\newenvironment{pf}{\medskip\noindent{\bf Proof.}\enspace}
{\hfill\newline\smallskip}

\begin{document}

\title	{\large\bf {A non-standard numerical scheme for an age-of-infection epidemic model}}

\author{
	E. Messina, \\
	{\small Department of Mathematics and Applications,}\\ {\small
		University of Naples ``Federico II''}\\
	{\small  Via Cintia, I-80126 Napoli, Italy, Member of the INdAM Research group 
		GNCS }\\
	{\small eleonora.messina@unina.it}\\ \\
	M. Pezzella\\
		{\small Department of Mathematics and Applications,}\\ {\small
		University of Naples ``Federico II''}\\
	{\small  Via Cintia, I-80126 Napoli, Italy, Member of the INdAM Research group 
		GNCS }\\
	{\small mario.pezzella@unina.it}	\\ \\
	A. Vecchio\\
	{ \small C.N.R. National Research Council of Italy,}\\ {\small Institute for 
		Computational Application ``Mauro Picone''}\\{\small Via P. Castellino,111 - 
		80131 Napoli -
		Italy, Member of the INdAM Research group GNCS}\\{\small 
		antonia.vecchio@cnr.it}}

\date{}
\maketitle

\begin{abstract}
We propose a numerical  method for approximating integro-dif\-fe\-ren\-tial equations arising in age-of-infection epidemic models. The method is based on a non-standard finite differences approximation of the integral term appearing in the equation. The study of convergence properties and the analysis of the qualitative behavior of the numerical solution show that it preserves all the basic properties of the continuous model with no restrictive conditions on the step-length $h$ of integration and that it recovers the continuous dynamic as  $h$ tends to zero. 
\end{abstract}

\vspace{.5cm}



\pagebreak

\section{Introduction}
\label{sec:intro}
Mathematical models based on non-linear integral and integro-differential equations are gaining increasing attention in mathematical epidemiology  due to their ability to incorporate the past infection dynamic into its current development \cite{ArinoChap2006,BrauerArt2017a,BrauerArt2017b,BrauerChap2016,BredaArt2012,LegacyBook,FengArt2000}. This property is particularly suitable to represent the evolution  of  diseases  where 
the dependence of infectivity on the time since becoming infected
plays a crucial role.
These \emph{age-of-infection} models, which contain an integral term describing the contribution of infected individuals to the total infectivity, need, in general, the support of numerical simulations for a complete qualitative  understanding and quantitative description. For this reason attention should be paid to set up  a numerical framework  that allows to provide real-time and reliable answers.  
In epidemic models governed by non-linear ordinary differential equations, numerical methods based on non-standard discretizations  are widely used, because they respond well to both requests. With respect to traditional methods, which may fail to capture some of the essential qualitative features of the model for certain values of the step-size, they allow a reliable description for different parameter values of the model. We  extend this approach to the integro-differential equation representing the  Kermack and McKendrick age-of-infection epidemic model: 
\begin{equation}\label{AgeOfInfectionModelIntegroDiff}
	S^{\prime}(t)=\beta S(t)\int_{0}^{\infty}[S^{\prime}(t-s)]A(s)\ ds,
\end{equation}
for which we refer to \cite{BrauerBook}, and the bibliography therein.
Here, $S(t)$ is the number of susceptibles at time $t$ and the constant $\beta$ represents the rate of effective contacts. Furthermore, $A(s)=\pi(s)B(s)$ is the  mean infectivity of members of the population with infection age $s,$ where $B(s)$ is the fraction of infected members remaining infected at infection age $s,$ and $0\leq\pi(s)\leq 1$  is the mean infectivity at infection age $s.$ \\
 Classical numerical approaches, like direct quadrature or collocation methods \cite{BrunnerBook2004,bvdh86},
give accurate approximations to the solution of equation \eqref{AgeOfInfectionModelIntegroDiff}, for sufficiently small values of the stepsize. However, there are two main crucial points, related to the nature of problem \eqref{AgeOfInfectionModelIntegroDiff}, that we want to underline here: the small stepsize needed to recover the continuous dynamic might be 
too demanding in terms of computational cost, furthermore
a result concerning the asymptotic behavior  of the numerical model, which parallels the one of the continuous problem, may be difficult to obtain.
\\ 
Our aim is to draw up a non-standard numerical scheme which
 preserves positivity and, in general, the dynamics of the continuous model \eqref{AgeOfInfectionModelIntegroDiff}. We show that
the non-standard method is dynamically consistent with the original continuous-time model and, therefore, it is expected to be a robust and efficient tool to integrate  problems with  more complex dynamics. This paper is organized as follows: in Section \ref{sec:model} we report the main results on the age-of-infection model  \eqref{AgeOfInfectionModelIntegroDiff} as developed in \cite{BrauerBook}. Then, in Section \ref{sec:nsfd} we formulate the numerical model, we give consistency and convergence results, and prove that it preserves  properties like positivity, monotonicity, boundedness,  for any value of the stepsize $h.$ The asymptotic dynamic of the numerical model is analysed in Section \ref{sec:dynamics}, where we propose discrete equivalents for the parameters characterizing the epidemics for fixed $h>0,$ and show the convergence to their continuous counterparts as $h\to 0.$ Finally, numerical experiments are reported in Section \ref{sec:num}, to show the theoretical results obtained, and some remarks in Section \ref{sec:concluding} conclude the paper.

\section{The age-of-infection epidemic model}
\label{sec:model}
We consider the Kermack and McKendrick age-of-infection epidemic model \eqref{AgeOfInfectionModelIntegroDiff}.
\\
It is assumed that the disease outbreak begins at time $t=0,$ so that $S(t)=N,$ for $t<0,$ and that there are no disease deaths, so that the total population size is a constant $N$ (see \cite{BrauerArt2008}).
If we introduce the total infectivity $\varphi(t)=-\int_{0}^{+\infty}{A(s)S^{\prime}(t-s)ds},$ at time $t,$ equation  \eqref{AgeOfInfectionModelIntegroDiff} can be rewritten as follows 
\begin{equation}\label{AgeOfInfectionModel}
\begin{split}
S^{\prime}(t)&=-\beta S(t)\varphi(t)\\
\varphi(t)&=\varphi_0(t)+ \beta \int_{0}^{t}A(t-s) S(s)\varphi(s) \ ds,
\end{split}
\end{equation}
where the function $\varphi_0(t)=-\int_{t}^{\infty}A(s)S^{\prime}(t-s) \ ds$ is the total infectivity, at time $t,$ of members of the population who were infected at  $t=0.$\\
From now on, we will refer equivalently to one or the other form of the model as needed, and we will base our investigation starting from the assumption that (see for example \cite{BrauerArt2017a})
\begin{equation}\label{eq:ipA}
	\int_0^{+\infty}A(s)ds<+\infty,
\end{equation}
and that (see for example \cite{BrauerArt2008}) all initial infectives have infection-age zero at $t=0,$
\begin{equation}\label{eq:ipPhi0}
	\varphi_0(t)=(N-S_0)A(t),
\end{equation}
where $S_0=S(0).$
There is ample literature, see  \cite{BrauerBook}, which deals with the description and analysis of age-of-infection epidemic models of the form \eqref{AgeOfInfectionModelIntegroDiff}. Here we outline the main facts, that will represent our guidelines for constructing a dynamics-preserving numerical scheme.
\begin{itemize}
	\item $S$ is a non-negative,  non-increasing function of time, and  decreases to a limit $S_{\infty}>0.$
	\item $\int_0^{\infty}\varphi(t)dt<+\infty,$ and  $\varphi(t)\to 0,$ for $t\to+\infty.$
	\item The \textit{basic reproduction number},
	\begin{equation}\label{R0}
		R_0=\beta N  \int_{0}^{\infty}A(s) \ ds,
	\end{equation}
is the number of secondary disease cases produced by one typical primary case, and represents an important indicator of the risk of epidemic; its role can be clarified through the following \textit{invasion criterion}.\\
	Based on the consideration that at the disease outbreak the entire population  is susceptible ($S(t)\approx N$), the linearization around the disease free equilibrium yields,
	\begin{equation*}
		S^{\prime}(t)=\beta N\int_{0}^{\infty}[S^{\prime}(t-s)]A(s)\ ds,
	\end{equation*}
	which has a solution $S(t)=Ne^{rt},$ with an exponential growth rate $r,$ if
	\begin{equation*}
		1=\beta N \int_{0}^{\infty} A(s)e^{-rs} \ ds.
	\end{equation*}
	Therefore $R_0$ can be expressed in terms of the initial growth rate
	\begin{equation}\label{relazione_r_esponenziale}
		R_0=\dfrac{\int_{0}^{\infty}A(s) \ ds}{ \int_{0}^{\infty} A(s)e^{-rs} \ ds},
	\end{equation}
	and it points out that an epidemic situation, for which initially the solution grows exponentially (see \cite{BrauerChap2016}), is characterized by
	\begin{equation*}
		r>0 \; \Leftrightarrow \; R_0>1.
	\end{equation*}
\item The \textit{final size} relation for the epidemic is
\begin{equation}\label{FinalSizeRelation}
	\log{\dfrac{S_0}{S_\infty} }=R_0\left(1-\dfrac{S_\infty}{N} \right),
\end{equation}
it has a unique solution  and gives the fraction $S_{\infty}/N$ of the population that escapes the epidemic.
\end{itemize}
The basic reproduction number and the final size relation represent important indicators for analyzing the behavior of epidemic models, therefore special attention will be given, in Section \ref{sec:dynamics}, to the analysis of their discrete counterparts.

\section{The non-standard finite difference scheme}
\label{sec:nsfd}
Consider an uniform mesh $t_n=nh,$ where $n=0,1,\dots,$ and $h>0$ is the stepsize.  
 We define the following discretization scheme for \eqref{AgeOfInfectionModel}
\begin{equation}\label{NSFDscheme}
	\begin{split}
	S_{n+1}=&S_n-h \beta S_{n+1}\varphi_n \\
	\varphi_{n+1}=&\varphi_0(t_{n+1})+h \beta\sum_{j=0}^{n} A(t_{n+1-j})S_{j+1}\varphi_j,
	\end{split}
	\end{equation}
for $n=0,1,\ldots,$ where $S_0=S(0),$  $\varphi_0=\varphi_0(0),$ and  $S_n\approx S(t_n),$  $\varphi_n\approx \varphi(t_n)$.
 Here, we have approximated the integral in \eqref{AgeOfInfectionModelIntegroDiff} by a  modified rectangular  rule, which is a left approximation in $A(t)\varphi(t)$ and a right approximation  in $S(t)$. For this reason, the numerical scheme  \eqref{NSFDscheme} falls into the class of non-standard finite difference methods, originally introduced for differential equations (see \cite{MickensBook}, and references therein) and only recently extended to integral problems \cite{LubumaArt}.

\subsection{Convergence}
\label{sec:convergence}
In this section we refer, when needed, to the equivalent compact notations for the continuous  problem and the numerical method,  respectively, 
\begin{equation}\label{VectorizedModel}
\begin{split}
	\begin{bmatrix}
	S(t) \\ \varphi(t)
	\end{bmatrix}&=\begin{bmatrix}
	S_0 \\ \varphi_0(t)
	\end{bmatrix}+ \beta \bigintss_{0}^{t} \begin{bmatrix}
	-1 & 0\\0 & A(t-s)
	\end{bmatrix} \begin{bmatrix}
	S(s)\varphi(s) \\ S(s)\varphi(s)
	\end{bmatrix} \ ds,
\end{split}
\end{equation}
and
\begin{equation}\label{VectorizedScheme}
\begin{split}
	\begin{bmatrix}
	S_{n+1} \\ \varphi_{n+1}
	\end{bmatrix}&=\begin{bmatrix}
	S_0 \\ \varphi_0(t_{n+1})
	\end{bmatrix}+ h \beta  \sum_{j=0}^{n}\begin{bmatrix}
	-1 & 0\\0 & A(t_{n+1}-t_j)
	\end{bmatrix} \begin{bmatrix}
	S_{j+1}\varphi_j \\ S_{j+1}\varphi_j
	\end{bmatrix} ,
\end{split}
\end{equation}
$n=0,1,\ldots.$ The analysis of the local error,
	\begin{equation}\label{LocalTruncationError}
	\begin{split}
		\delta(h;t_{n})=&\bigintss_{0}^{t_{n}} \begin{bmatrix}
			-1 & 0\\0 & A(t_{n}-s)
		\end{bmatrix} \begin{bmatrix}
			S(s)\varphi(s) \\ S(s)\varphi(s)
		\end{bmatrix} \ ds\\& -h \sum_{j=0}^{n-1} \begin{bmatrix}
			-1 & 0\\0 & A(t_{n}-t_j)
		\end{bmatrix} \begin{bmatrix}
			S(t_{j+1})\varphi(t_j) \\ S(t_{j+1})\varphi(t_j)
		\end{bmatrix},
	\end{split}
\end{equation}
 is not straightforward due to the non-standard nature of the integration rule. 
So, we need to prove the following  result.

\begin{lemma}\label{LemmaNSFDconsistent}
	Assume that the given function $A(t),$ describing problem \eqref{AgeOfInfectionModel}, is continuously differentiable on an interval $[0,T]$, with $T<+\infty,$ then the  method  \eqref{NSFDscheme} is consistent with \eqref{AgeOfInfectionModel}, of order $1$.
\end{lemma}
\begin{pf}
	The assumption on $A(t)$ implies that also $S(t)$ and $\varphi(t)$ are continuously differentiable on $[0,T].$ Let $h=T/M,$ with $M$ positive integer,  because of the convergence properties of rectangular quadrature rules (see for example \cite{DavisRabBook}), for each $j=0,\ldots, M-1$ it is
	\begin{equation}\label{eq:quaderr}\left\|
		\bigintss_{t_j}^{t_{j+1}} \begin{bmatrix}
			-S(s+h)\varphi(s) \\A(t_{n}-s)S(s+h)\varphi(s)
		\end{bmatrix}  ds-h \begin{bmatrix}
			-S(t_{j+1})\varphi(t_j) \\A(t_{n-j})S(t_{j+1})\varphi(t_j)
		\end{bmatrix}\right\|\leq ch^2,
	\end{equation}
	where $c>0$ does not depend on $h.$
For $n=0,\ldots,M,$ by simple manipulations, we write the integral in
 \eqref{LocalTruncationError} as:
	\begin{equation*}
		\begin{split}
		&\sum_{j=0}^{n-1}\bigintss_{t_j}^{t_{j+1}} \begin{bmatrix}
		-S(s)\varphi(s) \\A(t_{n}-s)S(s)\varphi(s)
		\end{bmatrix} ds\\
		=&\sum_{j=0}^{n-1}\bigintss_{t_j}^{t_{j+1}}\left(  \begin{bmatrix}
		-S(s+h)\varphi(s) \\A(t_{n}-s)S(s+h)\varphi(s)
		\end{bmatrix}
		-h \begin{bmatrix}
		-S^{\prime}(s+\theta_j h) \varphi(s) \\ A(t_{n}-s)S^{\prime}(s+\theta_j h)\varphi(s)
		\end{bmatrix} \right)  ds,
		\end{split}
	\end{equation*}
	with $0<\theta_j<1.$ 
	As a consequence, due to \eqref{eq:quaderr} and also to the regularity of the given functions, for the local truncation error \eqref{LocalTruncationError}, the bound
	\begin{equation*}
		\max_{0\leq n\leq M}\|\delta(h;t_{n})\|\leq Ch,
	\end{equation*}
holds. 
The positive constant $C$ depends on the bounds in $[0,T]$ for the functions and derivatives involved,  as well as on $T,$ but  not on $h.$ Thus, the proof is completed.
\end{pf}
Denote  by $e(h;t_n)=		\begin{bmatrix}
	S(t_{n}) \\ \varphi(t_{n})
\end{bmatrix}-\begin{bmatrix}
	S_{n} \\ \varphi_{n}
\end{bmatrix}$ the global error of the discretization \eqref{VectorizedScheme}. The following theorem, that can be easily proved by standard numerical techniques, provides sufficient conditions for the convergence of the numerical method.
\begin{theorem}\label{THMConvergence}
	Assume that the given function $A(t),$ describing problem \eqref{AgeOfInfectionModel}, is continuously differentiable on an interval $[0,T]$, and that $S_n,$ $\varphi_n$ are the approximations to \eqref{AgeOfInfectionModel},  defined by \eqref{NSFDscheme}, then
	\[\lim_{h\to 0}\max_{0\leq n\leq M}\|e(h;t_n)\|=0.\]
Furthermore, the order of convergence is $1.$\\
\end{theorem}
\subsection{Basic properties}
\label{sec:basics}
Consider the following result, whose proof comes immediately from  \cite[Lem.1]{MV2017}.
\begin{lemma}\label{th:lem1}
Let \begin{equation}
	\label{eq:sAp}
	A^{\prime}(t)\in L^1[0,+\infty),
\end{equation}
then the  quadrature error \begin{equation}\label{eq:tau} \tau(h)=\int_0^{+\infty}A(t)dt-h\sum_{n=0}^{+\infty}A(t_{n+1}),
\end{equation} tends to zero as $h\to 0.$ 
\end{lemma}
Observe that this lemma represents a generalization to the convergence result stated in \cite[Cor. p.208]{DavisRabBook}.\\
	From now on we assume that the given function $A(t)$  describing problem \eqref{AgeOfInfectionModel}, satisfies 
 \eqref{eq:sAp}. This also  implies that (see \cite[Lem.1]{MV2017})
\begin{equation}\label{eq:Abound}
		h\sum_{n=0}^{+\infty}A(t_{n+1})\leq \int_0^{+\infty}A(t)dt+h\bar A,
\end{equation}
with $\bar A=\int_0^{+\infty}|A^{\prime}(t)|dt.$
\begin{theorem}\label{THMmimimcNSFD}
 Let $(S_n,\;\varphi_n)$ be the solution to  the discrete equation \eqref{NSFDscheme}, with $h\geq 0,$  and non-negative initial values $S_0$ and $\varphi_0=\varphi_0(0)$.\\ 
Then:
	\begin{enumerate}
		\itemsep 0em
		\item $S_n$ and $\varphi_n$ are non-negative, $\forall n=1,2,\ldots,$
		\item the sequence $\{S_n\}_{n \in \mathbb{N}_0}$ is  non-increasing,
		\item  $\{S_n\}_{n \in \mathbb{N}_0}$ and $\{\varphi_n\}_{n \in \mathbb{N}_0}$ are bounded sequences, 
		\[\lim_{n\to \infty}S_n=S_{\infty}(h)\geq 0,\]
		and
		 \[\lim_{n\to \infty}\varphi_n= 0.\]
	\end{enumerate}
\end{theorem}
\begin{pf}
	For items $1.$ and $2.,$ we proceed by induction to prove that the statement $S_{n+1}\geq0$, $\varphi_{n+1}\geq 0$ and $S_{n+1}\leq S_n$, holds for all $n \in \mathbb{N}_0$ and $h\geq 0.$ The case $n=0$ is true because the initial values are non-negative.
	Assume that the properties are true for $j=1,\ldots,n-1,$ then:
	\begin{equation*}
		\begin{split}
			& 1+h \beta \varphi_{n} \geq 1 \ \Rightarrow \ 0\leq S_{n+1}=\dfrac{S_n}{1+h \beta \varphi_n}\leq S_n, \\
			& \varphi_{n+1}=\varphi_{0}(t_{n+1})+h \beta \sum_{j=0}^{n} A(t_{n+1-j})S_{j+1}\varphi_j\geq 0.
		\end{split}
	\end{equation*}
In order to prove item $3.$ observe that, for each $h\geq 0,$
since $\{S_n\}_{n \in \mathbb{N}_0}$ is a non-negative,  non-increasing sequence, then it is bounded from above by $S_0,$ and convergent to a finite  non-negative value. Furthermore, 
 the second of \eqref{NSFDscheme} and assumption \eqref{eq:ipPhi0} on $\varphi_0(t),$ imply that
$\displaystyle
\varphi_{n}\leq N\cdot\sup_{t\in [0,+\infty)}A(t),$ $n=0,1,\ldots,$
so also  $\{\varphi_n\}_{n \in \mathbb{N}_0}$ is bounded, by a constant $\displaystyle\bar{\varphi}=N\cdot\sup_{t\in [0,+\infty)}A(t)$ that does not depend on $h.$\\
Again from the second of \eqref{NSFDscheme} and assumption \eqref{eq:ipPhi0} we have
\begin{equation}\label{eq:sumfi}
	h\sum_{n=0}^{+\infty}\varphi_{n+1}=(N-S_0)h\sum_{n=0}^{+\infty}A(t_{n+1})+h \beta\sum_{j=0}^{+\infty} S_{j+1}\varphi_j h\sum_{n=0}^{+\infty} A(t_{n+1}).
\end{equation}
In equation \eqref{eq:sumfi}, the first of \eqref{NSFDscheme} and condition \eqref{eq:Abound} lead to
\[
h\sum_{n=0}^{+\infty}\varphi_{n+1}\leq (N-S_{\infty}(h))\left(\int_0^{+\infty}A(t)dt+h\bar A\right).
\]
It is clear that for $h$ greater than an arbitrary $\bar h>0,$ it is
\begin{equation}
	\label{eq:phibound}h\sum_{n=0}^{+\infty}\varphi_{n+1}\leq N\left(\int_0^{+\infty}A(t)dt+\bar h\bar A\right)<+\infty.
\end{equation}
Then $\varphi_n$  converges to zero, as $n\to +\infty,$ for any $h>0.$
\end{pf}
 Thus the properties of $S(t)$ and $\varphi(t),$ highlighted in Section \ref{sec:model}, are preserved by the numerical solution without any restriction on the stepsize $h.$
\section{Discrete asymptotic dynamics}
\label{sec:dynamics}
Since in this section we are going to study how  the numerical model \eqref{NSFDscheme} preserves the asymptotic dynamics of the continuous problem \eqref{AgeOfInfectionModel},  it is important to prove that the local error $\delta(h;t_n),$ given in \eqref{LocalTruncationError} is bounded  for  any $n=0,1,\ldots,$ and tends to zero as $h\to 0.$\\ 
We base our investigation  on the assumption \eqref{eq:sAp} on $A(t),$ which implies that Lemma \ref{th:lem1} and \eqref{eq:Abound} hold, and assures sufficient regularity for the solution to \eqref{AgeOfInfectionModel}.
\begin{theorem}\label{TEOerrLOCinf}
Assume that the given function $A(t),$ describing problem \eqref{AgeOfInfectionModel}, is continuously differentiable on $[0,+\infty).$ 
	 Then the scheme \eqref{NSFDscheme}, is consistent with  \eqref{AgeOfInfectionModel} on $[0,+\infty).$
\end{theorem}
\begin{pf}
	Consider
	\begin{equation*}
		\begin{split}
			&\bigintss_{t_{j-1}}^{t_j}\begin{bmatrix}
				-S(s+h)\varphi(s) \\ A(t_{n}-s)S(s+h)\varphi(s)
			\end{bmatrix} ds
		 -h \begin{bmatrix}
			-S(t_{j+1})\varphi(t_j) \\ A(t_{n}-t_j)S(t_{j+1})\varphi(t_j) \end{bmatrix}=\\		
	&	=\bigintss_{t_{j-1}}^{t_j}\begin{bmatrix}
				-S(s+h)\varphi(s)+S(t_{j+1})\varphi(t_j)  \\ A(t_{n}-s)S(s+h)\varphi(s)-A(t_{n}-t_j)S(t_{j+1})\varphi(t_j) 
			\end{bmatrix} ds\\
			=& -\bigintss_{t_{j-1}}^{t_j}\bigintss_{s}^{t_j} \dfrac{d}{d x}\begin{bmatrix}
				-S(x+h)\varphi(x) \\ A(t_{n}-x)S(x+h)\varphi(x)
			\end{bmatrix} dx \ ds .\\
		\end{split}
	\end{equation*}
Thus, for $n=0,1,\ldots,$
	\[\begin{split}
	&\left\|\bigintss_{0}^{t_{n}}\begin{bmatrix}
			-S(s+h)\varphi(s) \\ A(t_{n}-s)S(s+h)\varphi(s)
		\end{bmatrix} ds
		-h \sum_{j=0}^{n-1}\begin{bmatrix}
			-S(t_{j+1})\varphi(t_j) \\ A(t_{n}-t_j)S(t_{j+1})\varphi(t_j) \end{bmatrix}\right\|\\&\leq h \bigintss_{0}^{+\infty} \left\|\dfrac{d}{d x}\begin{bmatrix}
			-S(x+h)\varphi(x) \\ A(t_{n}-x)S(x+h)\varphi(x)
		\end{bmatrix} \right\|dx\\&+h\left\| \begin{bmatrix}
		-S(t_{n+1})\varphi(t_{n})+S(h)\varphi(0)\\ A(0)S(t_{n+1})\varphi(t_{n})-A(t_{n})S(h)\varphi(0)
	\end{bmatrix} \right\|\leq ch, \end{split}\]
where the constant $c>0$ does not depend on $n$ and $h,$ since, due to the assumption \eqref{eq:sAp} on  $A^{\prime}(t),$ it is $S^{\prime}(t),\;\varphi^{\prime}(t)\in L^1[0,+\infty),$ and $\varphi(t)$  bounded. Then, for the local error defined in \eqref{LocalTruncationError}, proceeding as in the proof of  Lemma \ref{LemmaNSFDconsistent}, we have
\[\sup_{n\geq 0}\|\delta(h;t_{n})\|\leq h\left(c+\sup_{n\geq 0}\bigintss_{0}^{t_{n}}\left\|\begin{bmatrix}
	S^{\prime}(s+\theta_nh)\varphi(s) \\ A(t_{n}-s)S^{\prime}(s+\theta_nh)\varphi(s)
\end{bmatrix} \right\|ds\right),\]
with $\theta_n\in (0,1),$  for all $n=0,1,\ldots.$ Since, $S^{\prime}(t)\in L^1[0,+\infty),$  $A(t)$ and $\varphi(t)$ are non-negative and bounded, the proof is completed. 
\end{pf}
We define 
\begin{equation}\label{eq:dr0}
	R_0(h)=h \beta N\sum_{n=0}^{+\infty}A(t_{n+1}),
\end{equation}
to be the numerical discretization of the basic reproduction number $R_0$ in \eqref{R0}. \\
Accordingly to the meaning of $R_0,$ the discrete reproduction number represents a threshold parameter for the numerical model. Indeed, the direct discretization to \eqref{AgeOfInfectionModelIntegroDiff} 
\begin{equation}\label{eq:disinv}
	\frac{S_{n+1}-S_{n}}{h}=h \beta S_{n+1}\sum_{j=0}^{\infty}A(t_{j+1})\frac{S_{n-j}-S_{n-j-1}}{h},
	\end{equation}
 is equivalent to \eqref{NSFDscheme}, with $\varphi_n=-h\sum_{j=0}^{+\infty}A(t_{j+1})\frac{S_{n-j}-S_{n-j-1}}{h},$ and \begin{equation}\label{eq:appophi0}
 	\varphi_0(t_n)=-h\sum_{j=n}^{+\infty}A(t_{j+1})\frac{S_{n-j}-S_{n-j-1}}{h}.
 \end{equation} Here, in agreement with the assumption, in Section \ref{sec:model}, that all initial infectives have infection age zero at $t=0,$  when computing $\varphi_0(t)$ at time $t_n,$ $A(t_j)=A(t_n),$ whenever $j>n.$ Thus, since for $n<0,$ $S_n=S(t_n)=N,$ the right-hand-side of \eqref{eq:appophi0} gives back expression \eqref{eq:ipPhi0} for $\varphi_0(t).$ \\If we assume that, initially, $S_n\approx N,$ the linearization of \eqref{eq:disinv} is
 \begin{equation}\label{eq:disinv_lin}
 	\frac{S_{n+1}-S_{n}}{h}=h \beta N\sum_{j=0}^{\infty}A(t_{j+1})\frac{S_{n-j}-S_{n-j-1}}{h},
 \end{equation}
which has an exponential solution $S_n=(1+r h)^n,$ if
\begin{equation}\label{eq:dinvasion}
	1=h \beta N\sum_{n=0}^{+\infty}A(t_{n+1})(1+rh)^{-(n+1)},\end{equation}
is satisfied. This can be interpreted as a discrete version of the \textit{invasion criterion}, where the solution to \eqref{NSFDscheme} initially grows exponentially as $(1+rh)^n,$ if $r>0.$ Furthermore, \eqref{eq:dr0} and \eqref{eq:dinvasion} give
\[R_0(h)=\frac{h\sum_{n=0}^{+\infty}A(t_{n+1})}{h\sum_{n=0}^{+\infty}A(t_{n+1})(1+rh)^{-(n+1)}},\]
which is the discrete equivalent to \eqref{relazione_r_esponenziale}, and
for which $r>0$ if and only if $R_0(h)>1.$ The fact that $R_0(h)=R_0-\beta N\tau(h),$ where $\tau(h)$ is the error defined in \eqref{eq:tau}, implies that the discrete scheme replicate the continuous dynamic for $h$ sufficiently small. \\
As already pointed out in Section \ref{sec:basics}, since $\{S_n\}_{n\in \mathbb N_0}$ is a non-negative, monotone, non-increasing sequence, it has a non-negative limit $S_{\infty}(h).$ From the first of \eqref{NSFDscheme} it is clear that
\begin{equation}\label{eq:Sinf}
	S_{\infty}(h)=\frac{S_0}{\prod_{n=0}^{\infty}(1+h \beta\varphi_n)},
\end{equation} thus implying the following relation for the discrete final size of the epidemic
\begin{equation}\label{eq:num_final1}\log{\frac{S_0}{S_{\infty}(h)}}=\sum_{n=0}^{\infty}\log{(1+h \beta\varphi_n)}.
\end{equation} 
The series at the right-hand side converges if and only if $h \beta\sum_{n=0}^{\infty}\varphi_n$ is finite. This is true because of \eqref{eq:phibound}.
Thus, the positiveness of $S_{\infty}(h)$ is guaranteed for any fixed value of the stepsize $h>0.$ Furthermore, we can express this series in terms of the numerical basic reproduction number \eqref{eq:dr0}, as follows:

\begin{theorem}\label{TEOfinalSIZE1}
Let $(S_n,\varphi_n)$ be the numerical solution to \eqref{AgeOfInfectionModel}, obtained  by the discrete scheme \eqref{NSFDscheme}, and
 $S_{\infty}(h)$ be defined in \eqref{eq:Sinf}. Then, for each $h>0,$ it holds:
		\begin{equation}\label{eq:fsize}
			h \beta \sum_{n=0}^{\infty}\varphi_{n+1}=(N-S_{\infty}(h))h \beta \sum_{n=0}^{+\infty}A(t_{n+1})=R_0(h)\left(1-\dfrac{S_\infty(h)}{N} \right) .
		\end{equation}
\end{theorem}
\begin{pf}
From \eqref{eq:sumfi}	and the first of \eqref{NSFDscheme} it is 
\[h \beta\sum_{n=0}^{+\infty}\varphi_{n+1}=\beta (N-S_0)h\sum_{n=0}^{+\infty}A(t_{n+1})-\beta (S_{\infty}(h)-S_0) h\sum_{n=0}^{+\infty} A(t_{n+1}).\]
Thus 
\begin{equation}\label{eq:24bis}
	h \beta\sum_{n=0}^{+\infty}\varphi_{n+1}=\beta (N-S_{\infty}(h))h\sum_{n=0}^{+\infty}A(t_{n+1}).
\end{equation}
From the definition of $R_0(h)$ in \eqref{eq:dr0}, \eqref{eq:fsize} is completely proved.
\end{pf}
From \eqref{eq:num_final1} and \eqref{eq:fsize}, we have
\[\log{\frac{S_0}{S_{\infty}(h)}}=U(h)R_0(h)\left(1-\dfrac{S_\infty(h)}{N}\right)+\mathcal O(h),\]
which is the discrete equivalent to the final size relation \eqref{FinalSizeRelation}, for any $h>0.$ The equivalence is more evident as $h\to 0,$ since the spurious term 
\[U(h)=\frac{\sum_{n=0}^{\infty}\log{(1+h \beta\varphi_n)}}{h \beta\sum_{n=0}^{\infty}\varphi_n},\] tends to $1,$ as shown in the next theorem.
	\begin{theorem}\label{VecchioProp1} Consider the solution $\varphi_n$ to the discrete equation \eqref{NSFDscheme}, obtained by a fixed stepsize $h>0,$ and define
		\begin{equation}\label{eq:u}
			u_n(h)=\dfrac{\log (1+h \beta \varphi_n)}{h \beta \varphi_n}.
		\end{equation}
	Then 
	\begin{description}
		\item[1.] $\lim_{h\to 0}u_n(h)=1,$ uniformly with respect to $n;$
		\item[2.] $\lim\limits_{h \to 0} \sum_{n=0}^{\infty} \left(\log (1+h \beta \varphi_n)-h \beta \varphi_n \right) =0.$
		\end{description}
	\end{theorem}
\begin{pf}
	From Theorem \ref{THMmimimcNSFD}, $\{\varphi_n\}_{n\in \mathbb N_0}$ is bounded by a constant $\bar \varphi,$ that does not depend on $h.$ 
	\\
Consider $h<(\beta \bar{\varphi})^{-1},$	by Taylor expansion in \eqref{eq:u} it is
\begin{equation*}
	u_n(h)=:\sum_{j=0}^{\infty}\dfrac{(-1)^j (h \beta \varphi_n)^j}{j+1}.
\end{equation*}
Thus, 
	\begin{equation*}
			|u_n(h)-1| \leq \sum_{j=1}^{\infty}(h \beta \varphi_n)^j\leq \sum_{j=1}^{\infty}(h \beta \bar{\varphi})^j= \dfrac{h \beta \bar{\varphi}}{1-h \beta \bar{\varphi}}.
	\end{equation*}
	Since the last term, in the previous inequality, tends to $0$ as $h$ goes to $0,$ it follows that $\forall \varepsilon>0,$ there exists $h_\varepsilon$ such that $h<h_\varepsilon$ implies $\dfrac{h \beta\bar{\varphi}}{1-h \beta \bar{\varphi}}<\varepsilon$.
	Finally, if $\varepsilon>0$ is fixed, the choice for  $h<\min\left\lbrace h_\varepsilon,\dfrac{1}{\beta \bar{\varphi}}\right\rbrace $ leads to $|u_n(h)-1|<\varepsilon$ $\forall n \in \mathbb{N}$. This proves the first part of the theorem.
\\	
	So, for $h<\bar{h}_{\varepsilon}$ and $n \in \mathbb{N},$ it is $1-\varepsilon<u_n(h)<1+\varepsilon$ and, since
	\begin{equation*}
		\sum_{n=0}^{\infty} \log(1+h \beta \varphi_n)=h \beta\sum_{n=0}^{\infty} \varphi_nu_n(h)\ 
	\end{equation*}
with $h \beta\sum_{n=0}^{\infty} \varphi_n<\infty,$ as proved in \eqref{eq:phibound},
	we can state that:
	\begin{equation*}
		(1-\varepsilon)h \beta\sum_{n=0}^{\infty} \varphi_n<h \beta\sum_{n=0}^{\infty} \varphi_nu_n(h)=\sum_{n=0}^{\infty} \log(1+h \beta \varphi_n)<(1+\varepsilon)h \beta\sum_{n=0}^{\infty} \varphi_n,
	\end{equation*}
	\begin{equation*}
		\Rightarrow \; \lim\limits_{h\to 0} \left(\dfrac{\sum_{n=0}^{\infty} \log(1+h \beta \varphi_n)}{h \beta\sum_{n=0}^{\infty} \varphi_n} \right)=1,
	\end{equation*}
	which completes the proof since, from \eqref{eq:phibound}, the  denominator is bounded.\\
\end{pf}
 As $h\to 0,$ we expect to recover the continuous dynamic. 
We prove the following theorem.
\begin{theorem}\label{VecchioTHM} Let $(S_n,\varphi_n)$ be the numerical solution to \eqref{AgeOfInfectionModel}, obtained  by the discrete scheme \eqref{NSFDscheme}, and
	$S_{\infty}(h)$ be defined in \eqref{eq:Sinf}. It holds:
	\begin{equation*}
		\lim\limits_{h\to 0}\left(\log{\dfrac{S_0}{S_\infty(h)} }-h \beta(N-S_\infty(h))\sum_{n=0}^{\infty}A(t_{n+1}) \right)=0.
	\end{equation*}
\end{theorem}
\begin{pf}
	We use \eqref{eq:num_final1} and simple algebraic manipulations, to obtain
	\[\begin{split}&\log{\frac{S_0}{S_{\infty}(h)}}-\beta (N-S_{\infty}(h))h\sum_{n=0}^{\infty}A(t_{n+1})=\\
		&\left(\sum_{n=0}^{\infty}\log{\left(1+h \beta\varphi_n\right)}-h \beta\sum_{n=0}^{\infty}\varphi_n\right)\\
		&+\left(h \beta\sum_{n=0}^{\infty}\varphi_n-\beta (N-S_{\infty}(h))h\sum_{n=0}^{\infty}A(t_{n+1})\right).\end{split}\]
	The statement is proved by taking the limit for $h\to 0$ of both terms and by using the result of Theorem \ref{VecchioProp1} and \eqref{eq:24bis}.
\end{pf}
\begin{remark}
Assume that there exists 
$\lim_{h\to 0}S_{\infty}(h)=S_{\infty}^*>0.$ For $h\to 0,$ $S_{\infty}$ and $S_{\infty}^*$ satisfy the same final size relation \eqref{FinalSizeRelation}, which  (see \cite{BrauerBook} or \cite{BrauerArt2008}) has a unique solution in $[0,S_0].$ 
Thus implying, once again, that the dynamic of the continuous model \eqref{AgeOfInfectionModel} is preserved by the numerical one \eqref{NSFDscheme}, for $h$ sufficiently small.
\end{remark}

\section{Numerical Examples}
\label{sec:num}
\begin{figure}[htp]
	\begin{center}
		\includegraphics[width=5in]{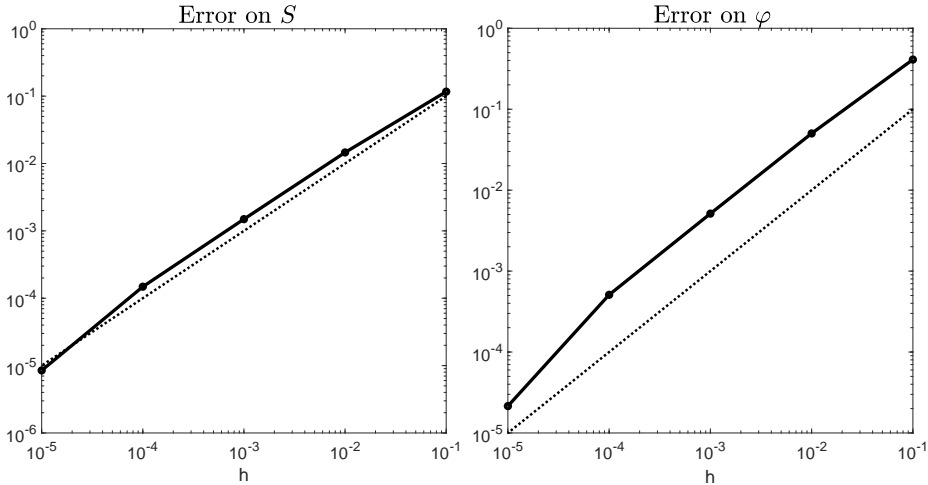}\\
		\caption{Problem \eqref{AgeOfInfectionModel}-\eqref{eq:test1}: norm of the relative errors (solid line) as functions of the stepsize,  compared to the slope of order one (dotted line).}\label{fig:Errors}
	\end{center}
\end{figure}

\begin{figure}[htp]
	\begin{center}
		\includegraphics[width=4in]{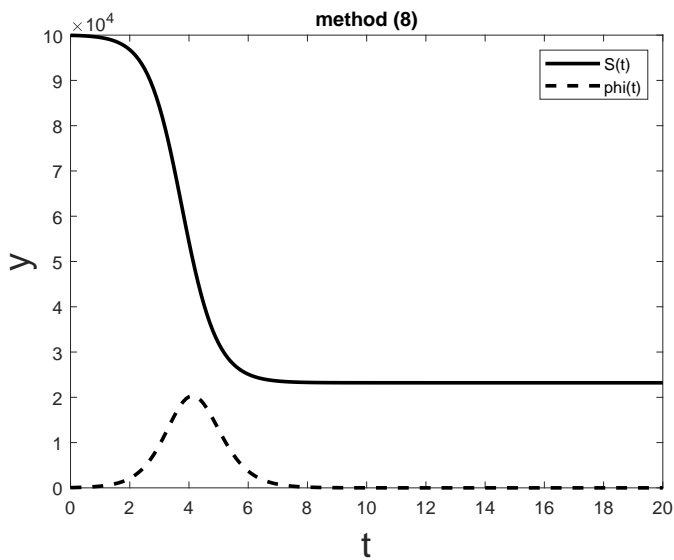}\\
		\caption{Problem \eqref{AgeOfInfectionModel}-\eqref{eq:test2}: numerical solution with $h=0.1.$ }\label{fig:test2}
	\end{center}
\end{figure}

In this section we report some numerical examples in order to show experimentally the theoretical results proved in the previous sections. 
For our experiments we choose illustrative test equations of the form \eqref{AgeOfInfectionModel}  and we use the non-standard method \eqref{NSFDscheme}.\\
As first example we integrate problem \eqref{AgeOfInfectionModel} for $t\in [0,1],$ with
\begin{equation}\label{eq:test1}
	A(t)=\dfrac{1}{(1+t)^2},\;\;\; N=10, \;\;\; S_0=9,  \;\;\;  \beta=0.3,
\end{equation}
and $\varphi_0(t)$ given by \eqref{eq:ipPhi0}.  
Theorem \ref{THMConvergence} states first order convergence of our scheme, 
and we observe, in Table \ref{tab:Errors and Experimental order} and Figure \ref{fig:Errors},  that the reduction of numerical errors as function of the stepsize confirms that behavior. Here, we have used the numerical solution computed with stepsize $h=10^{-6}$ as reference solution.
\begin{table}
	\centering
	\begin{tabular}{|c|c|c|c|c|}
		\hline
		$h$ & Error on $S$ & Error on $\varphi$ & Exp. ord. for $S$ & Exp. ord. for $\varphi$  \\
		\hline
		$10^{-1}$ & $1.17\cdot10^{-1}$ & $4.11\cdot10^{-1}$ & $\setminus\setminus $ & $\setminus\setminus$ \\
		$10^{-2}$ & $1.46\cdot10^{-2}$ & $5.02\cdot10^{-2}$ &$0.90$  & $0.91$ \\
		$10^{-3}$ & $1.49\cdot10^{-3}$ & $5.13\cdot10^{-3}$ &$0.99$  & $0.99$ \\
		$10^{-4}$ & $1.48\cdot10^{-4}$ & $5.09\cdot10^{-4}$ &$1.00$  & $1.00$ \\
		\hline
	\end{tabular}
	\bigskip
	\caption{Error values and experimental order of convergence for example \eqref{AgeOfInfectionModel}-\eqref{eq:test1}.} \label{tab:Errors and Experimental order}
\end{table}
In order to show the long time behavior of the numerical solution, we consider problem \eqref{AgeOfInfectionModel}, with a gaussian distribution for the total infectivity
\begin{equation}\label{eq:test2}
	A(t)=\frac{1}{\sigma\sqrt{2\pi} }e^{-\frac{(t-\mu)^2}{2\sigma^2}},\;\;\; N=100000, \;\;\; S_0=99950, \;\mu=0.2, \;\; \sigma=2\mu,
\end{equation}
we choose $\beta=3\cdot 10^{-5},$ and $\varphi_0(t)$ given by \eqref{eq:ipPhi0}. In Figure \ref{fig:test2} the behavior of the numerical solution is reported for $h=0.1.$  Here, it is clear that an epidemic occurs, according to the fact that the estimated value for the basic reproduction number, computed by \eqref{eq:dr0} is $R_0(h)\approx 1.29.$ 
By running the code on a sufficiently large interval, we obtain   the values reported in Table \ref{tab:test2} for $S_{\infty}(h)$
\begin{table}
	\centering
	\begin{tabular}{|c|c|}
		\hline
		$h$ & $S_{\infty}(h)$   \\
		\hline
		$10^{-1}$ & $2.3211\cdot10^{4}$ \\
		$10^{-2}$ & $1.8852\cdot10^{4}$ \\
		$10^{-3}$ & $1.8435\cdot10^{4}$ \\
		\hline
	\end{tabular}
	\bigskip
\caption{Values of the final size $S_{\infty}(h)$ as function of $h$ for problem \eqref{AgeOfInfectionModel}-\eqref{eq:test2}.} \label{tab:test2}
\end{table}
 which confirm the result in Theorem \ref{VecchioTHM}, compared to the value $S_{\infty}=1.8389\cdot10^{4},$ obtained by iteratively solving the non-linear final size relation \eqref{FinalSizeRelation}.\\
  \begin{figure}[htp]
 	\begin{center}
 		\includegraphics[width=5in]{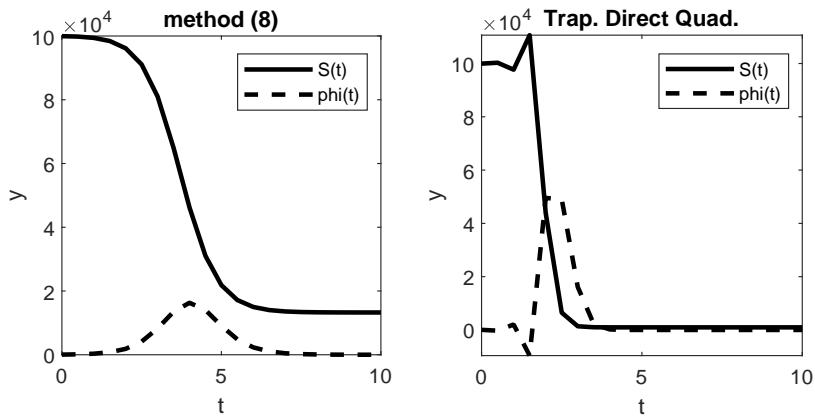}\\
 		\caption{Problem \eqref{AgeOfInfectionModel}-\eqref{eq:test2}: comparison of numerical solutions with $h=0.5.$ }\label{fig:test3}
 	\end{center}
 \end{figure}
 In our last experiment we compare the performances of the non-standard scheme \eqref{NSFDscheme} and a direct quadrature Trapezoidal method for problem \eqref{AgeOfInfectionModel}-\eqref{eq:test2} with $\beta=6\cdot 10^{-5},$ and using a relatively large stepsize $h=0.5.$ It is clear in Figure \ref{fig:test3} that the Trapezoidal Direct Quadrature method fails to preserve the positivity and monotonicity of the solution.

\section{Concluding remarks}
\label{sec:concluding}
In this work we study a numerical method for the integration of   age-of-infection epidemic models. Since these models have a great potential in the description of current epidemics, attention must be paid on the construction of discretization techniques that preserve the qualitative behavior of the continuous time model.  For the method we have proposed, that uses a non-standard discretization for the integral term characterizing the mathematical equations, we have conducted a comprehensive analysis which has allowed to ensure that the numerical solution is dynamically consistent with the continuous one, for any value of the step length. Furthermore, the method can be implemented in an explicit form and hence is computationally inexpensive.  The drawback, however stands in the convergence, which is linear,  and thus a severe restriction on the stepsize may be required if an accurate numerical answer is needed at finite time. 
Not wanting to give up the  robustness and the simplicity of implementation of these techniques, we can think of an adaptive strategy for the selection of the step-length in cases where the function $A$ rapidly grows and then becomes smoother after a certain period of time. This will be the subject of a future work, as well as a deeper study into the dynamic properties of classical higher order numerical approaches.
\section*{Acknowledgments} This work  was supported by~GNCS-INDAM.

\end{document}